# Optimal properties of some Bayesian inferences

**M. Evans and M. Shakhatreh**

*Dept. of Statistics, University of Toronto*
*e-mail:* mevans@utstat.utoronto.ca

**Abstract:** Relative surprise regions are shown to minimize, among Bayesian credible regions, the prior probability of covering a false value from the prior. Such regions are also shown to be unbiased in the sense that the prior probability of covering a false value is bounded above by the prior probability of covering the true value. Relative surprise regions are shown to maximize both the Bayes factor in favor of the region containing the true value and the relative belief ratio, among all credible regions with the same posterior content. Relative surprise regions emerge naturally when we consider equivalence classes of credible regions generated via reparameterizations.

## 1. Introduction

Suppose we have the ingredients for a proper Bayesian analysis. For this we observe data $x$ from a statistical model $\{f_\theta : \theta \in \Theta\}$, where $f_\theta$ is a density with respect to support measure $\mu$ on the sample space $\mathcal{X}$, and we have a proper prior density $\pi$ on $\theta$, with respect to support measure $\upsilon$ on $\Theta$. With these ingredients we have available the joint distribution of $(\theta, X)$, as given by the density $f_\theta(x)\pi(\theta)$ with respect to support measure $\nu \times \mu$, and the observed value $x$. We denote the prior predictive measure of $X$ by $M(B) = E_\Pi(P_\theta(B))$ and the posterior measure of $\theta$ by $\Pi(A \mid x)$. For a quantity of interest $\tau = \Upsilon(\theta)$, taking values in a set $\mathcal{T}$, we denote the marginal posterior and prior measures of $\Upsilon$ by $\Pi_\Upsilon(\cdot \mid x)$ and $\Pi_\Upsilon$ respectively, with corresponding densities $\pi_\Upsilon(\cdot \mid x)$ and $\pi_\Upsilon$, taken with respect to a support measure $\nu_\mathcal{T}$ on $\mathcal{T}$.

Bayes theorem, or the principle of conditional probability, says that any probability statements about the unknown $\theta$, after observing $x$, should be based on the posterior $\Pi(\cdot \mid x)$. These ingredients alone however, do not prescribe what $\gamma$-credible region $B_\gamma(x) \subset \mathcal{T}$ we should quote for $\tau = \Upsilon(\theta)$. Since there are typically many subsets of $\mathcal{T}$ containing $\gamma$ of the posterior probability, we need a rule for choosing among them.

Relative surprise credible regions for $\tau$, as discussed in Evans (1997), are based on a particular approach to assessing a hypothesis $H_0 : \tau = \tau_0$. For this we compute the *observed relative surprise* (ORS) given by

$$\Pi_\Upsilon \left( \frac{\pi_\Upsilon(\tau \mid x)}{\pi_\Upsilon(\tau)} > \frac{\pi_\Upsilon(\tau_0 \mid x)}{\pi_\Upsilon(\tau_0)} \,\bigg|\, x \right). \quad (1)$$

We see that (1) compares the relative increase in belief for $\tau_0$, from *a priori* to *a posteriori*, with this increase for each of the other possible values in $\mathcal{T}$. Other approaches to measuring surprise are discussed in Good (1988). For estimation, we consider (1) as a function of $\tau_0$ and select a value which minimizes this quantity as the estimate, called the *least relative surprise estimate* (LRSE). To obtain a $\gamma$-credible region for $\tau$ we simply invert (1) in the standard way to obtain the $\gamma$-*relative surprise region*

$$C_\gamma(x) = \left\{ \tau_0 \in \mathcal{T} : \Pi_\Upsilon \left( \frac{\pi_\Upsilon(\tau \mid x)}{\pi_\Upsilon(\tau)} > \frac{\pi_\Upsilon(\tau_0 \mid x)}{\pi_\Upsilon(\tau_0)} \,\bigg|\, x \right) \leq \gamma \right\}. \quad (2)$$

One virtue of relative surprise inferences is that they are invariant under reparameterizations.

In Evans, Guttman and Swartz (2006) it was shown that relative surprise inferences possess an optimal property in the class of Bayesian inferences. In that development, (2) was taken as the basic concept. In particular, if we consider the class of all $\gamma$-credible regions for $\tau = \Upsilon(\theta)$, then the $\gamma$-relative surprise region for $\tau$ has the smallest prior content among all $\gamma$-credible regions for this quantity. Hypothesis assessments and estimates are derived from relative surprise regions in a direct way and so also possess optimal properties. The LRSE is obtained by taking the region with $\gamma = 0$ and the ORS is obtained as $\inf\{\gamma : \tau_0 \in C_\gamma(x)\}$. In section 2 we show that this optimal property has a direct interpretation in terms of minimizing the prior probability of covering a false value and argue that this is an appropriate way to assess repeated sampling properties in contexts where we have a proper prior. In section 3 we prove that, for relative surprise regions, the prior probability of covering a false value is always bounded above by the prior probability of covering the true value and so such sets are, in a generalized sense, unbiased.

As discussed in Evans and Zou (2002) and Evans, Guttman and Swartz (2006), there is a close connection between relative surprise inferences and Bayes factors. In section 3 we establish some results that deepen this connection and show that relative surprise inferences lead to optimal results when interpreted in terms of Bayes factors. In particular, we prove that a $\gamma$-relative surprise region $C_\gamma(x)$ for $\tau = \Upsilon(\theta)$ always has a Bayes factor in favor of the region containing the true value bounded below by unity and, moreover, the Bayes factor is maximized among all $\gamma$-credible regions for $\tau$ by $C_\gamma(x)$. Further, we introduce the relative belief ratio as an alternative method for measuring change in belief from *a priori* to *a posteriori*, and show that this is also bounded below by unity for relative surprise regions and that such regions maximize this quantity as well. The lower bound can be seen as a natural consistency requirement on inferences in the sense that, it would be odd to report a $\gamma$-credible region for $\tau$ for which our belief in the set containing the true value declined from *a priori* to *a posteriori.* While a decline in belief from *a priori* to *a posteriori* makes sense for any particular subset of $\mathcal{T}$, it doesn't make sense for our best report for a set supposedly containing the true value, as the data are suggesting otherwise. Further, the optimality result indicates that we are making the best use of the data, from this point-of-view, when we choose to use relative surprise regions. In section 4 we show that in quite general circumstances relative surprise regions arise as a member of an equivalence class of credible regions under reparameterizations. Further, we show that choosing among these regions is equivalent to choosing the measure we use to construct an hpd-like credible region, where hpd stands for highest posterior density. We argue that the most natural choice of this measure is $\Pi_\Upsilon$, which gives relative surprise regions.

## 2. Covering false values

Suppose we have a rule for determining a $\gamma$-credible region for $\tau = \Upsilon(\theta)$ based on the sampling model and prior, i.e., for each $\gamma \in [0,1]$ and $x \in \mathcal{X}$, the rule determines a region $B_\gamma(x) \subset \mathcal{T}$ satisfying $\Pi_\Upsilon(B_\gamma(x) \,|\, x) \geq \gamma$. The coverage $P_\theta(\Upsilon(\theta) \in B_\gamma(X))$ of this region is then of considerable interest, particularly when $\Pi$ is taken to be a diffuse prior. In such a context it seems natural to ask that a $\gamma$-credible region satisfy, or at least approximately satisfy, the confidence property $P_\theta(\Upsilon(\theta) \in B_\gamma(X)) \geq \gamma$ for all $\theta \in \Theta$. In other words, in an i.i.d. sequence $x_i \sim P_\theta$ for $i = 1, 2, \ldots$, we require that the proportion of times that $\Upsilon(\theta) \in B_\gamma(x_i)$ is at least $\gamma$, and also that this property hold for all $\theta \in \Theta$. It is well-known that Bayesian credible regions do not generally possess this property, see Joshi (1974), and in fact can perform rather poorly in this regard. A number of papers discuss issues concerned with comparing frequency and Bayesian inferences including, Berger and Selke (1987), Casella and Berger (1987), and Samaniego and Reneau (1994) as well as the texts Gelman, Carlin, Stern and Rubin (2004), Carlin and Louis (2000), and Robert (2001).

We restrict to proper priors, as then, letting $E_M$ denote expectation with respect to the prior predictive distribution of the data,

$$\begin{aligned}
\gamma &\leq E_M(\Pi_\Upsilon(B_\gamma(X) \,|\, X)) = \int_\mathcal{X} \int_\Theta I_{B_\gamma(x)}(\Upsilon(\theta)) \, \Pi_\Upsilon(d\theta \,|\, x) \, M(dx) \\
&= \int_{\Theta \times \mathcal{X}} I_{B_\gamma(x)}(\Upsilon(\theta)) \, P_\theta(dx) \, \Pi(d\theta) = E_\Pi(P_\theta(\Upsilon(\theta) \in B_\gamma(X))).
\end{aligned} \quad (3)$$

This can be interpreted as saying that the prior probability the $\gamma$-credible region $B_\gamma$ contains a value $\Upsilon(\theta)$, when $\theta \sim \Pi$, is at least $\gamma$. This probability can also be given a long-run relative frequency interpretation in the i.i.d. sequence $(\theta_i, x_i) \sim \Pi \times P_\theta$ for $i = 1, 2, \ldots$, as the proportion of times $\Upsilon(\theta_i) \in B_\gamma(x_i)$. Various arguments can be offered for the restriction to proper priors, e.g., see DeGroot (1970). In particular, when we have a proper prior, this long-run relative frequency seems more appropriate than the confidence property, as the confidence property requires good coverage at values of $\theta$ that have *a priori* very little weight. Further, while the confidence property has its appeal, the plethora of absurd confidence regions, see Plante (1991) for some discussion, might at least lead one to doubt the wisdom of focusing too closely on confidence.

Property (3) holds for any $\gamma$-credible region $B_\gamma$ and so does not help us choose among them. Consider, however, the accuracy of the region $B_\gamma$, where this is measured by the probability of $B_\gamma$ covering an independent false value $\tau \sim \Lambda$ where $\Lambda$ is a probability measure on $\mathcal{T}$.

**Definition 1.** The *prior probability of covering a false value from probability measure* $\Lambda$ is given by

$$E_M(\Lambda(B_\gamma(X))) = \int_\Theta \int_\mathcal{T} P_\theta(\tau \in B_\gamma(X)) \, \Lambda(d\tau) \, \Pi(d\theta). \qquad (4)$$

So the "true" value of $\theta$ is generated from the prior $\Pi$, the data $x$ is generated from $P_\theta$, and the "false" value $\tau$ of the parameter of interest is generated from $\Lambda$ independent of the true value, i.e., $\tau$ has no connection with the data.

To obtain a $\gamma$-credible region $B_\gamma$ that minimizes (4) we make use of the following results. Suppose we have a probability measure $P$ and a $\sigma$-finite measure $Q$ on a set $\Omega$. Further, suppose that $P$ and $Q$ are both absolutely continuous with respect to the same measure on $\Omega$ with respective densities $p$ and $q$. Let $D_\gamma = \{\omega_0 \in \Omega : P(p(\omega)/q(\omega) > p(\omega_0)/q(\omega_0)) \leq \gamma\}$ for $\gamma \in [0,1]$. Lemma 1 and Theorem 2 are proved in Evans, Guttman and Swartz (2006). In that paper $P$ is taken to be the posterior but otherwise the proofs are the same.

**Lemma 1.** $P(D_\gamma) \geq \gamma$ with equality whenever the distribution of $p(\omega)/q(\omega)$, with $\omega \sim P$, has no atoms.

**Theorem 2.** *The set $D_\gamma$ minimizes $Q(D)$ among all measurable sets $D \subset \Omega$ satisfying $P(D) \geq P(D_\gamma)$. Further, when the distribution of $p(\cdot)/q(\cdot)$ has no atoms, then $D_\gamma$ minimizes $Q(D)$ among all measurable sets $D \subset \Omega$ satisfying $P(D) \geq \gamma$.*

The following result establishes the optimality, with respect to (4), of hpd-like credible regions as defined in (5).

**Theorem 3.** *Suppose that the probability distribution $\Lambda$ is also absolutely continuous with respect to $\nu_\mathcal{T}$ on $\mathcal{T}$ with density $\lambda$. Then, in the Bayesian model specified by $\Pi \times P_\theta$, the region $B_{\Lambda,\gamma}$ given by*

$$B_{\Lambda,\gamma}(x) = \left\{ \tau_0 \in \mathcal{T} : \Pi\left( \frac{\pi_\Upsilon(\tau \mid x)}{\lambda(\tau)} > \frac{\pi_\Upsilon(\tau_0 \mid x)}{\lambda(\tau_0)} \,\Big|\, x \right) \leq \gamma \right\} \qquad (5)$$

*minimizes (4) among all regions $B$ satisfying $\Pi_\Upsilon(B(x) \mid x) \geq \Pi_\Upsilon(B_{\Lambda,\gamma}(x) \mid x)$. If $\Pi_\Upsilon(B_{\Lambda,\gamma}(x) \mid x) = \gamma$ for each $x$, then $B_{\Lambda,\gamma}$ minimizes (4) among all $\gamma$-credible regions for $\Upsilon(\theta)$.*

*Proof.* Putting $P = \Pi_\Upsilon(\cdot \mid x)$ and $Q = \Lambda$ in Theorem 2, implies that $B_{\Lambda,\gamma}(x)$ minimizes $\Lambda(B(x))$ among all $B$ satisfying $\Pi_\Upsilon(B(x) \mid x) \geq \Pi_\Upsilon(B_{\Lambda,\gamma}(x) \mid x)$. Now $E_M(\Lambda(B(X)))$ is minimized, among all regions $B$ satisfying $\Pi_\Upsilon(B(x) \mid x) \geq$ $\Pi_\Upsilon(B_{\Lambda,\gamma}(x) \,|\, x)$ for each $x$, by $B = B_{\Lambda,\gamma}$. Observe that $E_M\left(\Lambda(B(X))\right) = \int_\Theta \int_\mathcal{X} E_\Lambda(I_{B(x)}(\tau)) \, P_\theta(dx) \, \Pi(d\theta) = \int_\Theta \int_\mathcal{X} \int_\mathcal{T} I_{B(x)}(\tau) \, \Lambda(d\tau) \, P_\theta(dx) \, \Pi(d\theta) = \int_\Theta \int_\mathcal{T} \int_\mathcal{X} I_{B(x)}(\tau) \, P_\theta(dx) \, \Lambda(d\tau) \, \Pi(d\theta) = \int_\Theta \int_\mathcal{T} P_\theta(\tau \in B(X)) \, \Lambda(d\tau) \, \Pi(d\theta)$ which is (4). This completes the proof. □

Specializing this to the choice $\Lambda = \Pi_\Upsilon$, we have the following result.

**Corollary 4.** *If $\Lambda = \Pi_\Upsilon$, then $B_{\Lambda,\gamma} = C_\gamma$, i.e., the optimal region is the $\gamma$-relative surprise region.*

This says that $C_\gamma$ minimizes the prior probability of covering a false value of the parameter of interest, when the false value follows the marginal prior distribution $\Pi_\Upsilon$ and is independent of the true value of the model parameter. It seems natural to take $\Lambda = \Pi_\Upsilon$ as this distribution identifies the values of $\tau$ that we consider *a priori* at least plausible. A repeated sampling interpretation of this is obtained by considering a sequence $(\theta_i, x_i, \tau_i)$, for $i = 1, 2, \ldots$, of independent values from the joint distribution $\Pi \times P_\theta \times \Pi_\Upsilon$. Then Corollary 4 says that, among all $\gamma$-credible regions $B_\gamma$ for $\Upsilon(\theta)$ formed from $\Pi \times P_\theta$, a $\gamma$-relative surprise region for $\Upsilon(\theta)$ minimizes the proportion of times the event $\tau_i \in B_\gamma(x_i)$ is true. Of course, the event $\Upsilon(\theta_i) \in B_\gamma(x_i)$ is also true at least $\gamma$ of the time in this sequence.

It is worth noting that Theorem 3 and Corollary 4 also hold when $\theta \sim \Pi^*$ for any probability distribution $\Pi^*$, i.e., $\int_\Theta \int_\mathcal{T} P_\theta(\tau \in B_\gamma(X)) \, \Lambda(d\tau) \, \Pi^*(d\theta)$ is minimized among all $\gamma$-credible regions $B_\gamma$ for $\Upsilon(\theta)$, constructed from $\Pi \times P_\theta$, by $B_{\Lambda,\gamma}$. The proof is the same. So, for example, we could take $\Pi^*$ to be degenerate at some value and $C_\gamma$ would still be optimal. From a practical, point-of-view, however, the choices $\Pi^* = \Pi$ and $\Lambda = \Pi_\Upsilon$ seem to be the most sensible, as (4) then has an interpretation as a prior probability.

In addition (4), with $\Lambda = \Pi_\Upsilon$, would appear to have several uses. First we can quote this probability as a way of assessing the accuracy of $C_\gamma$ with a given prior. If this probability is quite high, then we have a region with low accuracy. Also (4) can be used for experimental design purposes such as setting sample size. Consider the following example.

**Example 1** (Location normal). Suppose that $x = (x_1, \ldots, x_n)$ is a sample from the $N(\theta, 1)$ distribution and $\theta \sim N(0, \sigma^2)$ so the posterior distribution of $\theta$ is $N((n + 1/\sigma^2)^{-1} n\bar{x}, (n + 1/\sigma^2)^{-1})$. The ratio of the posterior density to the prior density is, in this case, proportional to the likelihood $\exp\left\{-n(\theta - \bar{x})^2/2\right\}$. Therefore, a $\gamma$-relative surprise interval for $\theta$ is a likelihood interval and so takes the form $C_\gamma(x) = \bar{x} \pm k_\gamma(n, \bar{x}, \sigma^2)$ where $k_\gamma(n, \bar{x}, \sigma^2) \geq 0$ satisfies

$$\begin{aligned}\gamma &= \Phi\big((n + 1/\sigma^2)^{-1/2} \bar{x}/\sigma^2 + (n + 1/\sigma^2)^{1/2} k_\gamma(n, \bar{x}, \sigma^2)\big) \\ &\quad - \Phi\big((n + 1/\sigma^2)^{-1/2} \bar{x}/\sigma^2 - (n + 1/\sigma^2)^{1/2} k_\gamma(n, \bar{x}, \sigma^2)\big). \end{aligned} \quad (6)$$

The value $k_\gamma(n, \bar{x}, \sigma^2)$ is easily obtained numerically from (6). Now $\Pi_\Upsilon(C_\gamma(x)) = \Phi\left((\bar{x} + k_\gamma(n, \bar{x}, \sigma^2))/\sigma\right) - \Phi\left((\bar{x} - k_\gamma(n, \bar{x}, \sigma^2))/\sigma\right)$ and $\bar{x} \sim N(0, \sigma^2 + 1/n)$. Therefore, when $\Lambda = \Pi_\Upsilon$, (4) is given by

$$\int_{-\infty}^{\infty} \Big[ \Phi\big(((\sigma^2 + 1/n)^{1/2}z + k_\gamma(n, (\sigma^2 + 1/n)^{1/2}z, \sigma^2))/\sigma\big)$$
$$- \Phi\big(((\sigma^2 + 1/n)^{1/2}z - k_\gamma(n, (\sigma^2 + 1/n)^{1/2}z, \sigma^2))/\sigma\big) \Big] \varphi(z) \, dz. \quad (7)$$

For example, the following table gives some values of the prior probability of covering a false value when $\gamma = .95$ and $\sigma^2 = 1$, based on a Monte Carlo integration sample size of $10^3$, with the standard errors in parentheses.

| $n$ | 1 | 10 | 25 | 50 |
| --- | --- | --- | --- | --- |
| $E_M\left(\Pi_\Upsilon(C_\gamma(X))\right)$ | .700 (.004) | .322 (.004) | .212 (.003) | .152 (.002) |

It is straightforward to show that (7) converges to 0 as $n \to \infty$. So, by choosing $n$ large enough, we can make (7) as small as we like and so control the error in our inference. If $\sigma^2 \to \infty$, so the prior is becoming more diffuse, the prior probability of covering a false value generated from the prior converges to 0. This is exactly how we would want our region to behave, namely, the data become much more important in determining the inference as the prior becomes more diffuse. For example, $C_\gamma(x) \to \bar{x} \pm n^{-1/2} z_{(1+\gamma)/2}$ as $\sigma^2 \to \infty$. So for a very diffuse prior, $C_\gamma(x)$ has a very small probability of covering an independently generated value from the prior.

Example 1 illustrates that we can't use (4), with $\Lambda = \Pi_\Upsilon$, to compare priors. A more concentrated prior will give a higher value for (4) than one more diffuse, however, we have different regions $C_\gamma(x)$ and different distributions for the false values under different priors. This emphasizes the importance of a careful choice of the prior so that unrealistic values of the parameter are excluded.

Similar optimality results can be obtained for the ORS given by (1). For suppose we agree to reject the hypothesis $H_0 : \Upsilon(\theta) = \tau_0$ whenever the ORS is greater than $\gamma$. This is equivalent to rejecting $H_0$ whenever $\tau_0 \in C_\gamma^c(x)$. Now consider the class of tests specified by $\gamma$-credible regions $B_\gamma$, so we reject $H_0$ whenever $\tau_0 \in B_\gamma^c(x)$. In this case, we want to find $B_\gamma$ maximizing

$$\int_\Theta \int_\mathcal{T} P_\theta(\tau \in B_\gamma^c(X)) \, \Pi_\Upsilon(d\tau) \, \Pi(d\theta \, | \, \Upsilon(\theta) = \tau_0). \quad (8)$$

This is the conditional prior probability, given that $H_0$ is true, that we would reject the hypothesis specified by $\tau$, when $\tau$ is a value independently generated from the prior. The quantity (8) is clearly analogous to power in the frequentist context. Then, arguing as in Theorem 3 and Corollary 4, we have that (8) is maximized by $C_\gamma^c(x)$ among all rejection regions with posterior content less than or equal to $1 - \gamma$. Also, we can use (8) to determine a sample size so that the test based on $C_\gamma^c$ has a prescribed value for this conditional prior probability.

## 3. Change in belief and unbiasedness

For $C \subset \mathcal{T}$, $BF_C(x) = \{\Pi_\Upsilon(C \, | \, x)/(1 - \Pi_\Upsilon(C \, | \, x))\}\{\Pi_\Upsilon(C)/(1 - \Pi_\Upsilon(C))\}^{-1}$ is the Bayes factor in favor of the true value of $\tau$ being in $C$. If we let $C$ shrink

nicely to $\tau_0$ (as in Rudin (1974), p. 163, a sequence of Borel sets $C_i$ shrinks nicely to a point $\tau_0$ if there is an $\alpha > 0$ such that each $C_i$ lies in an open ball $B(\tau_0, r_i)$ centered at $\tau_0$ and of radius $r_i > 0$, then $\mu(C_i) \geq \alpha\mu(B(\tau_0, r_i))$ for every $i$ where $\mu$ is volume measure, and $r_i \to 0$ as $i \to \infty$), then $BF_C(x)$ converges to $\pi_\Upsilon(\tau_0 \,|\, x)/\pi_\Upsilon(\tau_0)$ whenever these densities are continuous at $\tau_0$. So we can think of this quantity as an approximation to the Bayes factor associated with $\tau_0$ and the ORS is a calibration of this value to determine if it is indeed small and thus evidence against $\tau_0$ as a plausible value.

The Bayes factor in favor of $C$ is a measure of the change in our belief that $C$ contains the true value from *a priori* to *a posteriori*. Perhaps a simpler measure of this change in belief is given by the following.

**Definition 2.** The *relative belief ratio* of *a subset* $C \subset \mathcal{T}$, is given by $RB_C(x) = \Pi_\Upsilon(C \,|\, x)/\Pi_\Upsilon(C)$.

Again, as $C$ shrinks nicely to $\{\tau_0\}$, $RB_C(x)$ converges to $\pi_\Upsilon(\tau_0 \,|\, x)/\pi_\Upsilon(\tau_0)$. Note that $BF_C(x) = RB_C(x)/RB_{C^c}(x)$ and so $BF_C$ is not a function of $RB_C$ or conversely. They are measuring change in belief on different scales. Clearly the two will be approximately equal when $RB_{C^c}(x) \approx 1$ and this will occur whenever $C$ is "small".

Now consider a $\gamma$-relative surprise region $C_\gamma(x)$ for $\tau$. From (2), and the fact that the function $\Pi_\Upsilon(\pi_\Upsilon(\tau \,|\, x)/\pi_\Upsilon(\tau) > k \,|\, x)$ is right-continuous in $k$, there exists $k_\gamma(x)$ such that $C_\gamma(x) = \{\tau : \pi_\Upsilon(\tau \,|\, x)/\pi_\Upsilon(\tau) > k_\gamma(x)\}$. From this we have the following property for relative surprise regions. We assume throughout the remainder of this section that $\pi_\Upsilon(\tau) > 0$ for every $\tau \in \mathcal{T}$.

**Lemma 5.** *The relative surprise region* $C_\gamma(x)$ *satisfies* $RB_{C_\gamma(x)}(x) > k_\gamma(x)$.

*Proof.* We have that $\Pi_\Upsilon(C_\gamma(x) \,|\, x) = \int_{C_\gamma(x)} \pi_\Upsilon(\tau \,|\, x) \nu_\mathcal{T}(d\tau)$ which is clearly greater than $k_\gamma(x) \int_{C_\gamma(x)} \pi_\Upsilon(\tau) \nu_\mathcal{T}(d\tau) = k_\gamma(x) \Pi_\Upsilon(C_\gamma(x))$ proving the result. □

So Lemma 5 says that the ratio of posterior to prior probabilities of $C_\gamma(x)$ satisfies the same inequality that the respective densities do on this set.

We have an important lower bound on $BF_{C_\gamma(x)}(x)$ and $RB_{C_\gamma(x)}(x)$.

**Lemma 6.** *The relative surprise region* $C_\gamma(x)$ *satisfies* $BF_{C_\gamma(x)}(x) > 1$ *and* $RB_{C_\gamma(x)}(x) > 1$.

*Proof.* Clearly we have that $C_\gamma^c(x) = \{\tau : \pi_\Upsilon(\tau \,|\, x)/\pi_\Upsilon(\tau) \leq k_\gamma(x)\}$ and, as in Lemma 5, this implies that $RB_{C_\gamma^c(x)}(x) \leq k_\gamma(x)$. Combining this with Lemma 5 gives that $BF_{C_\gamma(x)}(x) > 1$. Since $BF_C(x) > 1$, then $1/\Pi_\Upsilon(C) > 1/\Pi_\Upsilon(C \,|\, x)$ and this implies that $RB_{C_\gamma(x)}(x) > 1$. □

Accordingly the Bayes factor and the relative belief ratio always indicate an increase in belief in the set $C_\gamma(x)$ from *a priori* to *a posteriori*. In particular, the posterior probability content of $C_\gamma(x)$ is always greater than its prior content. Note that, since $BF_{C^c}(x) = 1/BF_C(x)$, we have that $BF_{C_\gamma^c(x)}(x) < 1$ and $RB_{C_\gamma^c(x)}(x) < 1$ for a relative surprise region $C_\gamma(x)$.

Note that the fact the Bayes factor and relative belief ratio are always greater than 1 for a relative surprise region, does not imply that relative surprise inferences never find evidence against a hypothesized value $H_0 : \tau = \tau_0$. For we assess $H_0$ by computing (1), or equivalently from (2), computing $\gamma^* = \inf\{\gamma : \tau_0 \in C_\gamma(x)\}$. If $\gamma^*$ is large (near 1), then we have evidence against $H_0$. Alternatively, we could select an appropriate $\gamma$ and report $C_\gamma(x)$ as our best choice of a $\gamma$-credible region to contain the true value. If $\tau_0 \notin C_\gamma(x)$, then we have evidence against $H_0$.

Of course, there may be other credible regions with these properties. For example, hpd regions often have these properties, although there does not seem to be an easy general proof of this. In any case, the following shows that relative surprise regions are best from this point-of-view.

**Theorem 7.** *The set $C_\gamma(x)$ has maximal Bayes factor and maximal relative belief ratio among all measurable sets $C \subset \mathcal{T}$ satisfying $\Pi_\Upsilon(C \,|\, x) = \Pi_\Upsilon(C_\gamma(x) \,|\, x)$.*

*Proof.* From Theorem 2 we know that $\Pi_\Upsilon(C)$ is minimized, among all measurable $C$ satisfying $\Pi_\Upsilon(C \,|\, x) \geq \Pi_\Upsilon(C_\gamma(x) \,|\, x)$, by taking $C = C_\gamma(x)$. So $\Pi_\Upsilon(C)$ is also minimized by the same choice when we restrict to those $C$ satisfying $\Pi_\Upsilon(C \,|\, x) = \Pi_\Upsilon(C_\gamma(x) \,|\, x)$. Since $f(x) = (1-x)/x$ is decreasing in $x$, the result follows for the Bayes factor and is obvious for the relative belief ratio. □

Theorem 7 is most relevant when there are a number of credible regions, including the relative surprise region $C_\gamma(x)$, with posterior content exactly equal to $\gamma$. Theorem 7 then says that $C_\gamma(x)$ is the best choice among these regions from the point of view of the Bayes factor and the relative belief ratio, as it provides the largest increase in belief from *a priori* to *a posteriori*. We have the following immediate consequence.

**Corollary 8.** *Suppose that the true value of $\theta$ is selected according to $\Pi$. Then $E_M(BF_{B_\gamma(X)}(X))$ and $E_M(RB_{B_\gamma(X)}(X))$ are maximized, among credible regions $B_\gamma(x)$ satisfying $\Pi_\Upsilon(B_\gamma(x) \,|\, x) = \Pi_\Upsilon(C_\gamma(x) \,|\, x)$ for all $x$, by $B_\gamma(x) = C_\gamma(x)$.*

This says that the prior mean Bayes factor and prior mean relative belief ratio are maximized, by $C_\gamma(x)$.

Consider the following example as an illustration.

**Example 2** (Probability of joint success). Suppose we observe $x$ from a Binomial$(n, \theta_1)$, an independent $y$ from a Binomial$(n, \theta_2)$, we put independent uniform priors on $\theta_1$ and $\theta_2$ and we are interested in making inference about $\psi = \theta_1 \theta_2$. This is the probability of simultaneous success from tossing two coins where the coins have probability of heads equal to $\theta_1$ and $\theta_2$, respectively.

Suppose we have $n = 5$ and observe $x = 4$ and $y = 1$. In the following table we give some $\gamma$-hpd intervals and $\gamma$-relative surprise (rs) intervals for $\psi$. We see that these intervals are quite different. Also the relative surprise intervals always dominate the hpd intervals in the sense that the Bayes factor and relative belief ratio of the relative surprise interval are always greater than the corresponding quantities for the hpd interval, as proven generally in Theorem 7. The estimate determined by the hpd approach is the mode and this is given by .122 while the LRSE is .186. While the hpd intervals, in this example, always have RB $> 1$ and BF $> 1$, other methods of forming the intervals do not necessarily give intervals with these properties. For example, if we took the left-tail of the posterior as a $\gamma$-credible interval for $\psi$, then the left-tail .4-credible interval has RB $= .730$ and BF $= .640$.

| $\gamma$ | hpd | RB (hpd) | BF (hpd) | rs | RB (rs) | BF (rs) |
|---|---|---|---|---|---|---|
| .95 | (.008, .447) | 1.25 | 5.99 | (.028, .501) | 1.32 | 7.35 |
| .75 | (.032, .293) | 1.47 | 2.82 | (.071, .361) | 1.69 | 3.35 |
| .50 | (.059, .216) | 1.57 | 2.16 | (.110, .284) | 1.74 | 2.48 |
| .25 | (.089, .163) | 1.63 | 1.84 | (.119, .270) | 1.76 | 2.36 |

These results are somewhat typical for this context. For example, when $n = 20, x = 19, y = 19$, the .95-hpd interval is $(.675, .962)$ with RB $= 16.16$, and BF $= 305.40$, while the .95-relative surprise interval is $(.684, .990)$ with RB $= 16.98$, and BF $= 322.59$. In this case the posterior mode is .857 and the LRSE is .902.

In frequentist contexts, a confidence region is said to be unbiased, if the probability of the region containing a particular false value is always less than or equal to the probability of the region containing the true value. The following result shows that relative surprise regions are unbiased in a generalized sense.

**Theorem 9.** *For a relative surprise region, the prior probability of containing an independent value generated from the prior is always less than the prior probability of containing the true value, when it is generated from the prior.*

*Proof.* From Lemma 6 we have that $\Pi_\Upsilon(C_\gamma(x) \mid x) > \Pi_\Upsilon(C_\gamma(x))$ and so it follows that $E_M(\Pi_\Upsilon(C_\gamma(X))) < E_M(\Pi_\Upsilon(C_\gamma(X) \mid X)))$. By (4), $E_M(\Pi_\Upsilon(C_\gamma(X)))$ is the prior probability of $C_\gamma$ containing a false value while $E_M(\Pi_\Upsilon(C_\gamma(X) \mid X))) = E_\Pi(P_\theta(\Upsilon(\theta) \in C_\gamma(X)))$ is the prior probability that $C_\gamma(X)$ contains the true value $\Upsilon(\theta)$ when $\theta \sim \Pi, X \sim P_\theta$. □

## 4. Reparameterizations

A basic principle of inference is that inferences about a parameter of interest should be invariant under reparameterizations, e.g., whatever rule we use to obtain a $\gamma$-credible region $B_\gamma$ for a parameter of interest $\tau$, the rule should yield the region $\Psi B_\gamma$ for any 1-1, sufficiently smooth, reparameterization $\psi = \Psi(\tau)$. Relative surprise inferences satisfy this principle.

Suppose, however, that we insist on forming credible regions for parameters taking values in $\mathcal{T}$ by minimizing their $\Lambda$ content, where $\Lambda$ is also absolutely continuous with respect to $\nu_\mathcal{T}$ on $\mathcal{T}$ with density $\lambda$. Let $\mathcal{T}$ be an open subset of $R^k$ and $\mathcal{D}_{\mathcal{T},\mathcal{T}}$ denote the class of reparameterizations $\Psi : \mathcal{T} \to \mathcal{T}$ that are 1-1, onto, continuously differentiable and such that $\Psi^{-1}$ is continuously differentiable. Then, by Theorem 3, the $\gamma$-credible region for $\psi = \Psi(\tau)$ that has minimal $\Lambda$ content is given by the hpd-like region

$$B_{\Lambda,\gamma}^{\Psi}(x) = \left\{ \psi_0 \in \mathcal{T} : \Pi_\Upsilon \left( \frac{\frac{\pi_\Upsilon(\tau\,|\,x)J_\Psi^{-1}(\tau)}{\lambda(\Psi(\tau))}}{\frac{\pi_\Upsilon(\Psi^{-1}(\psi_0)\,|\,x)J_\Psi^{-1}(\Psi^{-1}(\psi_0))}{\lambda(\psi_0)}} > \;\Bigg|\; x \right) \leq \gamma \right\}, \quad (9)$$

where $J_\Psi(\tau)$ is the Jacobian of the transformation $\Psi$ evaluated at $\tau$.

Since $B_{\Lambda,\gamma}^{\Psi}(x)$ is a $\gamma$-credible region for $\psi = \Psi(\tau)$, then $\Psi^{-1}B_{\Lambda,\gamma}^{\Psi}(x)$ is a $\gamma$-credible region for $\tau$. Let $[B_{\Lambda,\gamma}(x)] = \{\Psi^{-1}B_{\Lambda,\gamma}^{\Psi}(x) : \Psi \in \mathcal{D}_{\mathcal{T},\mathcal{T}}\}$ be the class of $\gamma$-credible regions for $\tau$ that arise via reparameterizations, when using the measure $\Lambda$ to construct the credible regions. Each of the regions in $[B_{\Lambda,\gamma}(x)]$ is a plausible candidate as a $\gamma$-credible region for the parameter of interest and it is not clear how we should choose among them. The following result provides an approach to this choice.

**Lemma 10.** *For $\Psi \in \mathcal{D}_{\mathcal{T},\mathcal{T}}$ we have that $\Psi^{-1}B_{\Lambda,\gamma}^{\Psi}(x) = B_{\Lambda\circ\Psi,\gamma}(x)$ and so $[B_{\Lambda,\gamma}(x)] = \{B_{\Lambda\circ\Psi,\gamma}(x) : \Psi \in \mathcal{D}_{\mathcal{T},\mathcal{T}}\}$.*

*Proof.* For $A \subset \mathcal{T}$, we have that $\Lambda \circ \Psi(A) = \Lambda(\Psi(A)) = \int_{\Psi(A)} \lambda(\tau)\nu_\mathcal{T}(d\tau) = \int_A \lambda(\Psi(\tau))J_\Psi(\tau)\,\nu_\mathcal{T}(d\tau)$ and so the density of $\Lambda \circ \Psi$ is $\lambda(\Psi(\tau))J_\Psi(\tau)$. Therefore, by Theorem 3 and (9), the result follows. □

So it is equivalent to think of $[B_{\Lambda,\gamma}(x)]$ as containing all the $\gamma$-credible regions for $\tau$ obtained by minimizing the $\Lambda \circ \Psi$ content for some $\Psi \in \mathcal{D}_{\mathcal{T},\mathcal{T}}$. This result says that choosing among the elements of $[B_{\Lambda,\gamma}(x)]$ is equivalent to choosing which measure $\Lambda \circ \Psi$ we should use to optimize with respect to.

Now suppose that $\Lambda$ is a probability measure and define $\Psi_\Lambda : \mathcal{T} \to [0,1]^k$ so that $\Psi_\Lambda(\tau) \sim \text{Uniform}([0,1]^k)$ when $\tau \sim \Lambda$, e.g., we can take $\Psi_\Lambda$ to be the probability transform. Then, for probability measures $\Lambda_1$ and $\Lambda_2$, we can define $\Psi_{\Lambda_1,\Lambda_2} : \mathcal{T} \to \mathcal{T}$ by $\Psi_{\Lambda_1,\Lambda_2} = \Psi_{\Lambda_1} \circ \Psi_{\Lambda_2}^{-1}$ and thus $\Lambda_1 \circ \Psi_{\Lambda_1,\Lambda_2} = \Lambda_2$. Therefore, if $\Psi_{\Lambda_1,\Lambda_2} \in \mathcal{D}_{\mathcal{T},\mathcal{T}}$, we have that $[B_{\Lambda_1,\gamma}(x)] = [B_{\Lambda_2,\gamma}(x)]$. We say that $\Lambda_2$ is obtained via a smooth reparameterization from $\Lambda_1$ when $\Psi_{\Lambda_1,\Lambda_2} \in \mathcal{D}_{\mathcal{T},\mathcal{T}}$. We have the following result immediately from Corollary 4.

**Lemma 11.** *If $\Lambda$ is a probability measure and $\Psi_\Lambda^{-1} \circ \Psi_{\Pi_\Upsilon} \in \mathcal{D}_{\mathcal{T},\mathcal{T}}$, then $C_\gamma(x) \in [B_{\Lambda,\gamma}(x)]$.*

Note that, when $\Psi_{\Pi_\Upsilon}, \Psi_\Lambda$ are the respective probability transforms, $\lambda$ is continuous and positive and $\pi_\Upsilon$ is positive and continuous, then by the inverse function theorem, we must have that $\Psi_* = \Psi_\Lambda^{-1} \circ \Psi_{\Pi_\Upsilon} \in \mathcal{D}_{\mathcal{T},\mathcal{T}}$ and $J_{\Psi_*}(\tau) = \pi_\Upsilon(\tau)/\lambda(\Psi_*(\tau))$. Lemma 11 says that, in very general circumstances, when we choose to optimize with respect to a probability measure $\Lambda$ on $\mathcal{T}$, a relative surprise region is always available as an equivalent credible region under a reparameterization.

As previously noted, when we consider choosing among the elements of $[B_{\Lambda,\gamma}(x)]$ we need only consider which measure $\Lambda \circ \Psi$ is most appropriate. Theorem 3 says that choosing $\Lambda \circ \Psi$ leads to a region that minimizes the prior probability of covering a false value $\tau \sim \Lambda \circ \Psi$. Unless there are good reasons to do otherwise, the most appropriate weighting to apply to false values is given by the prior $\Pi_\Upsilon = \Lambda \circ \Psi_\Lambda^{-1} \circ \Psi_{\Pi_\Upsilon}$. This leads to a region that focuses on the parameter values that we believe are *a priori* important. For example, choosing a credible region that minimized the probability of covering false values that are well out of range where the prior placed most of its mass, would seem to be clearly inappropriate, as we presumably know *a priori* that these are unrealistic values.

The following illustrates the need for a rule to select a credible region.

**Example 3** (All $\gamma$-credible intervals can arise from reparameterizations). Suppose that the posterior of $\tau \in R^1$ is absolutely continuous, has finite second moment and $\pi_\Upsilon(\tau \,|\, x) > 0$ for every $\tau \in R^1$. Let $\tau_1 > \tau_0$ be such that $\Pi_\Upsilon((\tau_0, \tau_1) \,|\, x) = \gamma$. Let $\Lambda$ be a probability measure with density $\lambda(\tau) > 0$ for every $\tau \in R^1$. From Lemma 10, we have that the set of all $\Lambda$ hpd-like $\gamma$-credible intervals for $\tau$ obtained via reparameterizations, is the same as the set of all $\Lambda \circ \Psi$ hpd-like $\gamma$-credible intervals for $\tau$ as $\Psi$ ranges over all reparameterizations. Then, there is a constant $k_\gamma(x, \Psi)$ such that $B_{\Lambda \circ \Psi, \gamma}(x) = \{\tau : \pi_\Upsilon(\tau \,|\, x)/\lambda_\Psi(\tau) > k_\gamma(x, \Psi)\}$ where $\lambda_\Psi$ is the density of $\Lambda \circ \Psi$. Clearly, there exists a $\Psi$ such that $\lambda_\Psi(\tau) \propto \pi_\Upsilon(\tau \,|\, x)(1 + (\tau - (\tau_0 + \tau_1)/2)^2)$, i.e., the probability measure with this density is a smooth reparameterization of $\Lambda$, and so $B_{\Lambda \circ \Psi, \gamma}(x) = (\tau_0, \tau_1) \in [B_{\Lambda, \gamma}(x)]$.

While the reparameterization in the example depends on the data, it is not clear generally how to rule out reparameterizations. With the relative suprise rule this is not an issue, because of invariance.

So far we have restricted the discussion to probability measures $\Lambda$. Suppose, however, that $\Lambda$ is a bounded measure on $\mathcal{T}$. It is immediate, for any positive constant $b$, that $B_{b\Lambda, \gamma}(x) = B_{\Lambda, \gamma}(x)$. So we can take $b = 1/\Lambda(\mathcal{T})$ and simply treat $\Lambda$ as a probability measure, as we get the same set of credible regions and $C_\gamma(x) \in [B_{\Lambda, \gamma}(x)]$. Suppose now that $\Lambda$ is an unbounded measure with density $\lambda$ with respect to $\upsilon_\mathcal{T}$. Further suppose that there is a sequence of bounded measures $\Lambda_n$ with densities $\lambda_n$ with respect to $\upsilon_\mathcal{T}$, such that $\lambda_n \to \lambda$ pointwise as $n \to \infty$. For example, if $\Lambda$ and $\upsilon_\mathcal{T}$ are volume measure on $R^k$, then $\lambda \equiv 1$ and we can take $\lambda_n$ to be $(2\pi n)^{k/2}$ times a $N_k(0, nI)$ density. Then, when the posterior distribution of $\pi_\Upsilon(\tau \,|\, x)/\lambda(\tau)$ is continuous, we have that $B_{\Lambda_n, \gamma}(x) \to B_{\Lambda, \gamma}(x)$ as $n \to \infty$, since $\liminf B_{\Lambda_n, \gamma}(x) = \limsup B_{\Lambda_n, \gamma}(x) = B_{\Lambda, \gamma}(x)$ up to a set having posterior measure 0. If $\lambda$ and $\pi_\Upsilon$ are positive and continuous, then we have that $C_\gamma(x) \in [B_{\Lambda_n, \gamma}(x)]$ for each $n$. Therefore, $B_{\Lambda, \gamma}(x)$ is approximated by $B_{\Lambda_n, \gamma}(x)$ for large $n$ and $C_\gamma(x)$ is equivalent to this set under a reparameterization. Accordingly, we can think of $C_\gamma(x)$ as being approximately equivalent to $B_{\Lambda, \gamma}(x)$ under a reparameterization.

## 5. Conclusions

Relative surprise regions have been shown to minimize the prior probability of covering a false value from the prior. This prior probability can be seen to serve as a measure of accuracy of a credible region and can be used for design purposes. Further, relative surprise regions have optimal properties with respect to the Bayes factor and relative belief ratio of the region. Finally, we have shown that relative surprise regions arise very naturally when we consider choosing among equivalent credible regions based on reparameterizations.

The relevance of our results in a particular application depends on the prior. In our view this is no different than concerns about the relevance of our choice of a sampling model in a problem, i.e., if we make a poor choice, then any inferences drawn based on this model are at least suspect. Model checking methods, can increase our confidence, when the model passes, that our choice makes sense. Similarly, methods for checking for prior-data conflict, such as those discussed in Evans and Moshonov (2006, 2007), can increase our confidence that the prior we have chosen makes sense. When the model and prior pass such checks, then optimal inferences drawn from such ingredients have greater force. In particular, the repeated sampling interpretations based upon the prior, then seem much more appropriate to us than the common frequentist practice of looking for procedures that possess good properties uniformly over all values of the model parameter, i.e., even at values of the parameter that we believe *a priori* are not relevant.

## Acknowledgements

The authors thank a referee and Associate Editor for comments that led to improvements in the paper.